\def\timestamp{%
Time-stamp: <example-vanmill.tex: Wednesday 15-12-2004 at 17:01:46 (cet)>}
\def\stripname Time-stamp: <#1 #2>{#2}
\edef\filedate{\expandafter\stripname\timestamp}
\documentclass{amsart}
\usepackage{amsrefs}

\theoremstyle{plain}
\newtheorem{theorem}{Theorem}[section]
\newtheorem{lemma}[theorem]{Lemma}

\theoremstyle{remark}
\newtheorem{remark}[theorem]{Remark}

\newcommand{\closure}{\overline}
\newcommand{\prodless}[2]{\ensuremath{#1^{<#2}}}

\newcommand{\cont}{\mathfrak{c}}
\newcommand{\pee}{\mathfrak{p}}

\newcommand{\axiom}{\mathsf}
\newcommand{\CH}{\axiom{CH}}
\newcommand{\MA}{\axiom{MA}}

\newcommand{\C}{\mathsf{C}}    
\newcommand{\N}{\mathbb{N}}    
\newcommand{\Z}{\mathbb{Z}}    

\DeclareMathSymbol\restr{2}{AMSa}{"16}

\newcommand{\orpr}[2]{\langle#1,#2\rangle}

\begin{document}

\title{A note on an example by van Mill}

\author{K. P. Hart}
\address{Faculty of Electrical Engineering, Mathematics, and
         Computer Science\\
         TU Delft\\
         Postbus 5031\\
         2600~GA {} Delft\\
         the Netherlands}
\email{K.P.Hart@EWI.TUDelft.NL}
\urladdr{http://aw.twi.tudelft.nl/\~{}hart}

\author{G. J. Ridderbos}
\address{Faculty of Sciences\\
         Division of Mathematics\\
         Vrije Universiteit\\
         De Boelelaan 1081\textsuperscript{a}\\
         1081~HV {} Amsterdam\\
         the Netherlands\\
	   Telephone number: (0031) 20 5987677}
\email{gfridder@cs.vu.nl}

\begin{abstract} 
Adapting an earlier example by J. van~Mill, we prove that there exists
a zero-dimensional compact space of countable $\pi$-weight and uncountable
character that is homogeneous under $\MA+ \neg\CH$, but not under~$\CH$.
\end{abstract}

\date{\filedate}


\maketitle

\section{Introduction}

In \cite{vanMill2003} van Mill constructed a compact Hausdorff space
of countable $\pi$-weight and character~$\aleph_1$ with the curious
property that $\MA+\neg\CH$ implies it is homogeneous whereas 
$\CH$~implies that it is \emph{not}.

The space was obtained as a resolution of the Cantor set, where every point
was replaced by the uncountable torus~$T^{\omega_1}$.
Close inspection of the proof of homogeneity from $\MA+\neg\CH$ reveals that
one needs any possible replacement, $Y$, of~$T^{\omega_1}$ to have 
the following properties:
\begin{enumerate}
\item $Y$ is compact Hausdorff and homogeneous;
\item the weight and character of $Y$ are $\aleph_1$;\label{wY+chiY}
\item there are an autohomeomorphism~$\eta$ of~$Y$ and a point~$d$ such that
         the positive orbit $\{\eta^n(d):n\in\omega\}$ is dense; and
\item $Y$ is a retract of~$\gamma\N$, whenever $\gamma\N$ is a 
      compactification of~$\N$ with~$Y$ as its remainder.\label{retract}
\end{enumerate}

In Section~\ref{sec:construction} of this note we shall show that that these
four properties do indeed suffice for a proof of homogeneity from
$\MA+\neg\CH$ and in Section~\ref{sec:input} we exhibit a zero-dimensional
space with these properties.
This then will establish that there is a zero-dimensional compact Hausdorff
space of countable $\pi$-weight and character~$\aleph_1$ with the curious
property that $\MA+\neg\CH$ implies it is homogeneous whereas $\CH$~implies
that it is \emph{not}.
Indeed, the power $2^{\omega_1\times \Z}$ with the mapping~$\eta$,
defined by $\eta(x)(\alpha,n)=x(\alpha,n+1)$, is as required.

\section{The Construction}
\label{sec:construction}

To keep our presentation self-contained we give an alternative description
of van Mill's construction.
We fix a compact space~$Y$ as in the introduction, along with the map~$\eta$
and the point~$d$; for each $n\in\omega$ we write $d_n=\eta^n(d)$.

The underlying set of $X$ is the product $\C\times Y$, where $\C$~is the
Cantor set~$2^\omega$.
Before we define the topology on~$X$ we fix some notation.

Given $s\in\prodless2\omega$, so $s$ is a finite sequence of zeroes and ones,
we put
$$
[s]=\{x\in\C:s\subseteq x\}.
$$
The family $\bigl\{[s]:s\in\prodless2\omega\bigr\}$ is the canonical base 
for the
topology of~$\C$.
If $s\in\prodless2\omega$ and $x\in\C$ then $s*x$ denotes the concatenation
of~$s$ and~$x$.

Given $x\in \C$ and $n\in\omega$ we put $U_{x,n}=[x\restr n]$, the $n$th
basic neighbourhood of~$x$,
and $C_{x,n}=U_{x,n}\setminus U_{x,n+1}$. 
Note that $C_{x,n}$ is of the form $U_{y,n+1}$ for some suitably 
chosen $y\in\C$.

Let $\orpr{x}{y}$ be a point of~$X$.
The basic neighbourhoods of~$\orpr{x}{y}$ will be those of the
form 
$$
U_{x,n}\otimes W = \bigl(  \{x\}\times W \bigr) \cup
        \bigcup\{C_{x,m}\times Y : m\ge n, d_m\in W\},
$$
with $n\in\omega$ and where $W$ runs through the neighbourhoods of~$y$ in~$Y$.

We will use this description of~$X$ but we invite the interested reader
to verify that $X$ is indeed a resolution of~$\C$, where each~$x$ is resolved
into~$Y$ via the map $f_x:\C\setminus\{x\}\to Y$ defined by
$f_x(y)=d_n$ iff $y\in C_{x,n}$, see~\cites{Fedorcuk1968,Watson} for details
on resolutions.

The following lemmas are easily verified and left to the reader.

\begin{lemma}
If $U_{x,n}$ and $W$ are clopen then so is $U_{x,n}\otimes W$,
hence $X$~is zero-dimensional if $Y$ is.\qed
\end{lemma}
 
\begin{lemma}
The character of $\orpr{x}{y}$ in~$X$ is the same as the character
of~$y$ in~$Y$, hence $\chi(X)=\aleph_1$.\qed
\end{lemma}

\begin{lemma}
The family $\bigl\{[s]\times Y:s\in\prodless2\omega\bigr\}$ is a $\pi$-base 
for~$X$, hence $\pi(X)=\aleph_0$.\qed
\end{lemma}

The next lemma should be compared with \cite{Watson}*{Theorem~3.1.33}.

\begin{lemma}
The space $X$~is compact Hausdorff.  
\end{lemma}

\begin{proof}
Let $\mathcal{U}$ be a basic open cover of~$X$.
Much like in a standard proof that the product of two compact spaces is 
compact one proves that around every vertical line one can put an open strip 
that is covered by finitely many elements of~$\mathcal{U}$.
Let $z\in X$. 
In the case that there is a set $U_{x,n}\otimes W\in\mathcal{U}$ such that
$z\in C_{x,m}\times Y$ for some~$m$ we are done.
In the other case for every $y\in Y$ there are $W_y$ and $n_y$ such that
$y\in W_y$ and $U_{z,n_y}\otimes W_y\in\mathcal{U}$; finitely many of these,
indexed by the finite set~$F$ say, will cover~$\{z\}\times Y$.
Let $n=\max\{n_y:y\in F\}$, then these sets cover the strip $U_{z,n}\times Y$
as well.
\end{proof}

At this point we know that $\CH$ implies that $X$~is not homogeneous,
because under $\CH$ we have $\chi(X)=2^{\pi(X)}$.
Corollary~1.2 in \cite{vanMill2003} forbids this for homogeneous compacta.

We now begin to work toward a proof that $X$~is homogeneous if
$\MA+\neg\CH$ is assumed.
As a first step we show that points with the same second coordinate
are similar.

\begin{lemma}
For each $a\in C$ the map $T_a$ defined by $T_a(x,y)=\orpr{x+a}{y}$
is an autohomeomorphism of~$X$.  
\end{lemma}

\begin{proof} 
One easily verifies that
$T_a[U_{x,n}\otimes W]=U_{x+a,n}\otimes W$, and this suffices.  
\end{proof}

The hard work will be in establishing that points with the same first
coordinate are similar.
We begin by showing that the special clopen sets $[s]\times Y$ are all 
homeomorphic and we give canonical homeomorphisms between them.

\begin{lemma}\label{lemma.xi-st}
Let $s,t\in\prodless2\omega$, put $k=|t|-|s|$ and define 
$\xi_{s,t}:[s]\times Y\to[t]\times Y$ by 
$\xi_{s,t}(s*x,y)=\orpr{t*x}{\eta^k(y)}$.
Then $\xi_{s,t}$ is a homeomorphism.  
\end{lemma}

\begin{proof}
It is not hard to show that
$\xi_{s,t}[U_{s*x,n}\otimes W]=U_{t*x,n+k}\otimes\eta^k[W]$, which suffices.
\end{proof}

In this lemma $k$ may be positive or negative and we need both possibilities,
see the last three lines of this section.

For ease of notation we let $e$ be the point of~$\C$ with all coordinates
zero and we abbreviate $U_{e,n}$ and $C_{e,n}$ by $U_n$ and $C_n$ respectively.
We shall prove, assuming $\MA+\neg\CH$,
that for every homeomorphism $f:Y\to Y$ there is a homeomorphism 
$\bar f:X\to X$ such that $\bar f(e,y)=\orpr{e}{f(y)}$
for all~$y\in Y$.
This will complete the proof that $X$ is homogeneous under this assumption
because it shows that all points of the form $\orpr ey$ are similar.

For $n\in\omega$ let $x_n\in\C$ be the point in $C_n$ with all coordinates
zero except for the $n$th.
Let $E=\{\orpr{x_n}{d_n}:n\in\omega\}$; note that $E$~is discrete. The proof
of the following lemma is that of~\cite{vanMill2003}*{Lemma 5.1}.

\begin{lemma}
$\closure E=E\cup\bigl(\{e\}\times Y\bigr)$.\qed
\end{lemma}

Thus, $\closure E$ is a compactification of~$\N$, whose remainder is~$Y$.

Let $f$ be an autohomeomorphism of~$Y$, considered to be acting 
on~$\{e\}\times Y$.
Theorem~4.2 from~\cite{vanMill2003}, an extension of a result due to
Matveev~\cite{Matveev}, now implies that there is an autohomeomorphism~$F$
of~$\closure E$ that extends~$f$.
This requires assumptions (\ref{wY+chiY}) and (\ref{retract}) on~$Y$.
It suffices to know that the weight of~$Y$ is less than the cardinal~$\pee$;
the inequality $\aleph_1<\pee$ follows from $\MA+\neg\CH$.

Note that the action of~$F$ on~$E$ is given by a permutation~$\tau$ 
of~$\omega$.

We use $\tau$ to define our extension~$\bar f$:
use the maps~$\xi_{s,t}$ from Lemma~\ref{lemma.xi-st} to map $C_n\times Y$
onto~$C_{\tau(n)}\times Y$ for all~$n$.
The verification that $\bar f$ is a homeomorphism is as in~\cite{vanMill2003}

\section{The Input}
\label{sec:input}

We now show that the space $Y=2^{\omega_1\!\times \Z}$ provides suitable input
for the construction in the previous section.

It is clear that $Y$ is zero-dimensional, compact and Hausdorff and that
its weight and character are~$\aleph_1$. Since $Y$ is a product of second
countable compacta, condition (\ref{retract}) is also satisfied, 
by~\cite{vanMill2003}*{Corollary~4.5}.

We need to find an autohomeomorphism $\eta$ and a point~$d$ such that
$\{\eta^n(d):n\in\omega\}$ is dense in~$Y$.
The map $\eta : Y\to Y$ is defined co-ordinatewise as follows 
(with $\alpha\in\omega_1$ and $i\in \Z$),
$$
\eta(y)(\alpha,i) = y(\alpha, i+1).
$$
We may think of points of $Y$ as $\omega_1$ by $\Z$ matrices. The action of
$\eta$ on such a matrix consists of shifting every row one step downwards.  
For $n<\omega$, $[-n,n]$ is the set $\{-n,-n+1, \ldots , n-1 , n\}$. 

To make the point~$d$ we take a countable 
dense subset $Q$ of~$2^{\omega_1\!\times \Z}$ 
(cf.~\cite{Engelking}*{Theorem 2.3.15}) 
and we enumerate $Q\times\omega$ as $\{\orpr{q_k}{n_k}:k<\omega\}$.
We define the point~$d$ by concatenating the restrictions 
$q_{k}\restr \omega_1 \times [-n_k,n_k]$.
First write $N_k=\sum_{j<k}\bigl(2\cdot n_k + 1\bigr)$ for all~$k$ and 
then define, for each~$\alpha$ and each~$n$
$$
d(\alpha,n)=\begin{cases}
                                 0  & n<0\\
         q_k (\alpha, n- N_k - n_k) & N_k\le n<N_{k+1}.\\
            \end{cases}
$$
Next we verify that the point $d$ has a dense positive orbit under $\eta$. 
Observe that it follows by construction that for all $k<\omega$ we have
$$ 
\eta^{N_k + n_k}(d)\restr (\omega_1 \times[-n_k,n_k])=
q_k \restr (\omega_1 \times [-n_k,n_k]).
$$

Now let an arbitrary basic open subset $U$ of $2^{\omega_1\times \Z}$ be given by
a function $s: F\to 2$, where $F\subseteq \omega_1\! \times \Z$ is finite. Thus $U$ is
given by
$$
U = \{ y \in 2^{\omega_1\! \times \Z} : s\subseteq y \}.
$$
Since $F$ is finite, we may find $n$ such that $F\subseteq \omega_1\! \times [-n,n]$.
The set $Q$ was chosen dense in $2^{\omega_1\! \times \Z}$, so there is a $k < \omega$
with $q_k \in U$ and $n_k = n$. It follows that
$$
q_k \restr (\omega_1 \times [-n_k,n_k]) \supseteq s
$$
from which it follows that
$$
\eta^{N_k + n_k}(d)\restr(\omega_1 \times[-n_k,n_k])\supseteq s.
$$
This implies that $\eta^{N_k+n_k}(d)\in U$.

We find that the set $\{\eta^{n}(d) : n<\omega \}$ is dense in
$2^{\omega_1\!\times \Z}$, which means that we are done.

\begin{remark}
This argument shows that also $2^{\cont \times \Z}$ has a point with a 
dense orbit under the shift, that is, not only is the Cantor cube $2^\cont$
separable, it even has an autohomeomorphism~$\eta$ and a point~$d$ whose
positive orbit is dense.
\end{remark}


\end{document}